\documentclass[leqno, 12pts]{article}

\usepackage{graphicx}

\usepackage{stmaryrd}
\usepackage{amsmath,amssymb,amsfonts}
\usepackage{times}
\usepackage[T1]{fontenc}
\usepackage{graphicx}
\usepackage[english]{babel}
\usepackage{amsmath}
\usepackage{amsthm}
\usepackage{amssymb}
\usepackage{enumerate}
\usepackage{xspace}
\usepackage{euscript}
\usepackage{graphicx}
\usepackage{amscd}
\usepackage{epsfig}
\usepackage{tabularx}
\usepackage{color}
\usepackage{bm}
\usepackage{tikz} 

\newtheorem{theorem}{Theorem}[section]

\theoremstyle{definition}

\theoremstyle{remark}
\newtheorem{remark}[theorem]{Remark}

\numberwithin{equation}{section}

\begin{document}

\title{Parallel Multi-Step Contour Integral Methods for Nonlinear Eigenvalue Problems}



\author{Yingxia Xi \thanks{School of Mathematics and Statistics Nanjing University of Science and Technology, Nanjing, 210094, China. ({\tt xiyingxia@njust.edu.cn})} 
\and Jiguang Sun \thanks{Department of Mathematical Sciences, Michigan Technological University, Houghton, MI 49931, U.S.A. ({\tt  jiguangs@mtu.edu}).}}
\date{}

\maketitle

\begin{abstract}
We consider nonlinear eigenvalue problems to compute all eigenvalues in a bounded region on the complex plane. Based on domain decomposition and contour integrals, two robust and scalable parallel multi-step methods are proposed. The first method 1) uses the spectral indicator method to find eigenvalues and 2) calls a linear eigensolver to compute the associated eigenvectors. The second method 1) divides the region into subregions and uses the spectral indicator method to decide candidate regions that contain eigenvalues, 2) computes eigenvalues in each candidate subregion using Beyn's method; and 3) verifies each eigenvalue by substituting it back to the system and computes the smallest eigenvalue. Each step of the two methods is carried out in parallel. Both methods are robust, accurate, and does not require prior knowledge of the number and distribution of the eigenvalues in the region. Examples are presented to show the performance of the two methods. 
\end{abstract}

\section{Introduction}
Let $T: \Omega \to \mathbb C^{n, n}$ be a matrix-valued function where $\Omega \subset \mathbb C$ is bounded.
We consider the nonlinear eigenvalue problem of finding $\lambda \in \mathbb C$ and ${\boldsymbol x} \in \mathbb C^n$ such that
\begin{equation} \label{GLambdax}
T(\lambda) {\boldsymbol x}= {\boldsymbol 0}.
\end{equation}
Nonlinear eigenvalue problems arise from various applications in science and engineering \cite{Betcke2013}. Of particular interests are the problems resulting from the numerical discretization of nonlinear eigenvalue problems of partial differential equations, e.g.  band structures for dispersive media \cite{Xiao2021JOSAA}, scattering resonances \cite{DyatlovZworski2019,MaSun2023AML}, transmission eigenvalue problems \cite{Gong2022MC,SunZhou2016}. In many cases, one cannot assume much prior knowledge on the distribution and number of eigenvalues.

Nonlinear eigenvalue problems have been an important topic in the numerical linear algebra community \cite{Voss2004BIT,BaiSun2005,Betcke2013,GuttelFisseur2017}. Many methods such as iterative techniques and linearization algorithms have been investigated. Recently, contour integral methods have become popular \cite{Asakura2009, Beyn2012,Gong2022MC,BrennanEtal2023}. The main ingredient of these methods is the complex contour integral for holomorphic matrix functions. Among them, the spectral indicator method (recursive integral method) is very simple and easy to implement \cite{Huang2016JCP,Gong2022MC}. It can be used as a standalone eigensolver or a screening tool for the distribution of eigenvalues. In contrast, Beyn’s algorithm uses Keldysh’s Theorem to retrieve all the spectral information and works well provided some knowledge of the number and distribution of eigenvalues \cite{Beyn2012}. 

In this paper, we propose two parallel multi-step contour integral methods, {\bf pmCIMa} and {\bf pmCIMb}. The goal is to find all the eigenvalues (and the associated eigenfunctions) of $T(\cdot)$ in $\Omega$ without much prior knowledge. {\bf pmCIMa} has two steps: 1) uses the spectral indicator method to compute eigenvalues and 2) call a linear eigensolver to compute the associated eigenvectors. {\bf pmCIMb} consists of three steps: (1) detection of eigenvalues in subregions of $\Omega$ using the spectral indicator method \cite{Huang2016JCP,Huang2018NLA}, (2) computation of eigenvalues using Beyn's method in the selected subregions \cite{Beyn2012}, and (3) verification of eigenvalues using a standard linear solver. Each step of the two methods is implemented in parallel. We demonstrate the performance of both methods using two nonlinear eigenvalue problems. 


The proposed methods have a couple of merits. They do not require prior knowledge of the number and distribution of the eigenvalues. The robustness are achieved by carefully choosing the parameters and adding a validation step. The parallel nature of domain decomposition guarantees the high efficiency. The accuracy is obtained by using a reasonable number of quadrature points due to the exponential convergence of the trapezoidal rule. The rest of the paper is organized as follows. Section \ref{Pre} contains the preliminary for contour integral methods. In particular, we introduce the spectral indicator method and Beyn's method. In Section \ref{pmSIM}, we present the algorithms for the parallel multi-step contour integral methods and discuss various implementation issues. The proposed methods are tested by two examples in Section~\ref{NE}.  Finally, we draw some conclusions in Section~\ref{CF}.

\section{Preliminary}\label{Pre}
We introduce SIM and Beyn's method, both of which use the contour integral, a classic tool in complex analysis. 
We refer the readers to \cite{Huang2016JCP,Huang2018NLA,Beyn2012} for more details.

\subsection{Spectral Indicator Method}
Let $\Omega$ be bounded and connected such that $\partial \Omega$ is a simple closed curve.  If $T(z)$ is holomorphic on $\Omega$, then $T(z)^{-1}$ is meromorphic on $\Omega$ with poles being the eigenvalues of $T(z)$. 

Assuming $T(z)$ has no eigenvalues on $\partial \Omega$, define an operator $P \in \mathbb C^{n, n}$ by
\begin{equation}\label{P}
P=\dfrac{1}{2\pi i}\int_{\partial \Omega}T(z)^{-1}dz,
\end{equation}
which is a projection from $\mathbb C^n$ to the generalized eigenspace associated with all the eigenvalues of $T(\cdot)$ in $\Omega$.
If there are no eigenvalues of $T(\cdot)$ inside $\Omega$, then $P = 0$, and $P{\boldsymbol f} = {\boldsymbol 0}$ for all
${\boldsymbol f} \in \mathbb C^n$. Otherwise, $P{\boldsymbol f} \ne {\boldsymbol 0}$ with 
probability $1$ for a random vector ${\boldsymbol f}$. Hence $P{\boldsymbol f}$ can be used to decide if $\Omega$ contains eigenvalues or not. 

Define an indicator for $\Omega$ using a random $ {\boldsymbol f} \in \mathbb C^n$ as
\begin{equation}\label{indicatorRIM}
I_\Omega:=\left |P {\boldsymbol f}\right |.
\end{equation}
SIM computes $I_{\Omega}$ and checks if $\Omega$ contains eigenvalues of $T(\cdot)$ by comparing $I_{\Omega}$ and a threshold value $tol_{ind}$. If $I_{\Omega} \ge tol_{ind}$, it divides $\Omega$ into subregions and compute the indicators for these subregions. The procedure continues until the regions are smaller than the required precision $tol_{eps}$. The original SIM is proposed for linear eigenvalue problems but can be directly applied to nonlinear eigenvalue problems \cite{Gong2022MC}. The following algorithm is recursive. One can change it to a loop easily.
\begin{itemize}
\item[] {\bf SIM}$(\Omega, tol_{eps}, tol_{ind})$
\item[]{\bf Input:}  region of interest $\Omega$, precision $tol_{eps}$, threshold $tol_{ind}$.
\item[]{\bf Output:}  eigenvalues $\lambda$ inside $\Omega$.
\item[1.] Compute ${I_{\Omega}}$ using a random vector ${\boldsymbol f}$.
\item[2.] If $I_{\Omega} < tol_{ind}$, exit (no eigenvalues in $\Omega$).
\item[3.] Otherwise, compute the diameter $h$ of $\Omega$.
			\begin{itemize}
				\item[-] If $h  > tol_{eps} $, 
				partition $\Omega$ into subregions $\Omega_j, j=1, \ldots J$.
						\begin{itemize}
						\item[] for $j=1$ to $J$
						\item[] $\qquad${\bf SIM}$(\Omega_j, tol_{eps}, tol_{ind})$.
						\item[] end
						\end{itemize}
				\item[-] else, 
						\begin{itemize}
							\item[] set $\lambda$ to be the center of $\Omega$. 
							\item[] output $\lambda$ and exit.
						\end{itemize}
			\end{itemize}
	
\end{itemize}

The major task of SIM is to approximate the indicator $I_\Omega$ defined in \eqref{indicatorRIM} using some quadrature rule
\begin{equation}\label{XLXf}
P{\boldsymbol f} \approx  \dfrac{1}{2 \pi i} \sum_{j=1}^{N} \omega_j {\boldsymbol x}_j,
\end{equation}
where $\omega_j$'s are quadrature weights and ${\boldsymbol x}_j$'s are the solutions of the linear systems
\begin{equation}\label{linearsys}
T(z_j){\boldsymbol x}_j = {\boldsymbol f}, \quad j = 1, \ldots, N.
\end{equation}
The linear solver is usually problem-dependent and out of the scope of the current paper. 
The total number of the linear systems \eqref{linearsys} to solve by SIM is at most
\begin{equation}\label{complexityRIM}
 2N \lceil \log_2(h/tol_{eps})\rceil k,
\end{equation}
where $k$ is the (unknown) number of eigenvalues in $\Omega$, $N$ is the number of the quadrature points, $h$ is the diameter of $\Omega$, $tol_{eps}$ is
the required precision, and $\lceil \cdot \rceil$ denotes the least larger integer. In general, SIM needs to solve a large number of linear systems. 

The original version of SIM does not compute eigenvectors. However, it has a simple fix. If $\lambda$ is found to be an eigenvalue, then $T(\lambda)$ has an eigenvalue $0$. The associated eigenvectors for $\lambda$ are the eigenvectors associated to the eigenvalue $0$ for $T(\lambda)$. Typical eigensovlers such as Arnodi methods \cite{arpack} can be used to compute the eigenvectors of $T(\lambda)$.

\subsection{Beyn's Method}
We briefly introduce Beyn's method which computes all eigenvalues in a given region~\cite{Beyn2012}.
Assume that there exist $k \le n$ eigenvalues $\{\lambda_j\}_{j=1}^k$ inside $\Omega$ and no eigenvalues lie on $\partial \Omega$. 

Let $V \in\mathbb{C}^{n,k}, k \leq n,$ be a random full rank matrix and
\begin{align}\label{integrals1}
C_0&=\frac{1}{2\pi i}\int_{\partial \Omega}T(z)^{-1}Vdz \in \mathbb{C}^{n \times k}, \\ 
\label{integrals2} 
C_1&=\frac{1}{2\pi i}\int_{\partial \Omega} z T(z)^{-1} V dz \in \mathbb{C}^{n \times k}.
\end{align}
Let the singular value decomposition of $C_0$ be given by
$$C_0=V_0\Sigma_0W_0^H,$$
where $V_0\in \mathbb C^{n \times k}$, $\Sigma_0=diag(\sigma_1,\sigma_2,\cdots,\sigma_k)$, $W_0\in \mathbb C^{k \times n}$. Then, the matrix 
\begin{equation}\label{D}
D:=V_0^HC_1W_0\Sigma_0^{-1}\in\mathbb C^{k,k}
\end{equation}
is diagonalizable with eigenvalues $\lambda_1,\lambda_2,\cdots,\lambda_k$ and associated eigenvectors ${\boldsymbol s}_1, \ldots {\boldsymbol s}_k$.
The eigenvectors for $T(\cdot)$ are given by ${\boldsymbol v}_j = V_0 {\boldsymbol s}_j$.

Consequently, to compute the eigenvalues (and eigenvectors) of $T(\cdot)$ in $\Omega$, one first computes \eqref{integrals1} and \eqref{integrals2} to obtain $C_0$ and $C_1$, then performs the singular value decomposition for $C_0$, and finally compute the eigenvalues and eigenvectors of $D$.


The number of eigenvalues inside $\Omega$, i.e., $k$, is usually unknown. One would expect that there is a gap between larger singular values and smaller singular values of $C_0$.
However, this is not the case if there exist eigenvalues (both inside and outside $\Omega$) close to $\partial \Omega$. This makes the rank test challenging. Furthermore, if $k$ is larger than $n$, some adjustment is needed (see Remark 3.5 of \cite{Beyn2012}).

%
%

\section{Parallel multi-step Contour Integral Method}\label{pmSIM}
We now propose two parallel multi-step contour integral methods, {\bf pmCIMa} and {\bf pmCIMb}, to compute the eigenvalues of $T(\cdot)$ inside a bounded region $\Omega \subset \mathbb C$. {\bf pmCIMa} 1) uses SIM in parallel to compute eigenvalues in $\Omega$ to a given precision, and 
2) plugs these values into $T(\cdot)$ and uses a linear eigensolver to compute the eigenvectors. {\bf pmCIMb}  consists of three steps: 1) detection of eigenvalues in the subregions using SIM in parallel, 2) computation of eigenvalues and eigenvectors using Beyn's method for each subregion, and 3) verification of the eigenvalues using a linear eigensolver. 

For both methods, we need to divide $\Omega$ into subregions.
It is convenient to cover $\Omega$ with squares. However, the quadrature error of the trapezoid rule for a circle decays exponentially with an exponent depending on the product of the number of quadrature points and the minimal distance of the eigenvalues to the contour \cite{Beyn2012}. The take this advantage, we first cover $\Omega$ by squares and then use the disks $\Omega_i$'s circumscribing these squares. Since there is a gap between the square and its circumscribing disk, eigenvalues outside the square are discarded.
 In fact, the gap between a square
and its circumscribing circle becomes negligible when the size of the square is small enough. In the rest of the paper, we shall not distinguish the square and its circumscribing circle for convenience.

\subsection{pmCIMa}
It is difficult to decide the threshold value $tol_{ind}$ for the indicator defined in \eqref{indicatorRIM}. 
Instead, we propose a more effective indicator similar to the one in \cite{Huang2018NLA}. Assume $N$ is an even number and let $P {\boldsymbol f}|_{N}$ be the approximation of ${\boldsymbol f}$ using $N$ quadrature points. Define an indicator, still denoted by $I_\Omega$, as
\begin{equation}
I_\Omega = \frac{1}{\sqrt{n}}\left | \frac{P {\boldsymbol f}|_N}{P {\boldsymbol f}|_{N/2}} \right|,
\end{equation}
where $P {\boldsymbol f}|_{N/2}$ is the approximation of $P {\boldsymbol f}$ using $N/2$ quadrature points $\{{\boldsymbol x}_2, {\boldsymbol x}_4, \ldots, {\boldsymbol x}_N\}$. It is expected that $I_\Omega \approx 1$ if there exists eigenvalues in $\Omega$ and $I_\Omega \approx e^{-Cn/2}$ for some constant $C$ otherwise. In the implementation, we use $N=16$ and set $tol_{ind} = 0.1$, which is quite reliable.

In Step 1, we cover $\Omega$ with subregions $\{\Omega_i\}_{i=1}^I$ and compute $I_{\Omega_i}$ to determine if $\Omega_i$ contains eigenvalues. If $I_{\Omega_i} > tol_{ind}$, then $\Omega_i$ is subdivided into smaller regions and these regions are saved for the next around. The procedure continues until the size of the region is smaller than the precision $tol_{eps}$. Then the centers of these small regions are the approximate eigenvalues.

In Step 2, an approximate eigenvalue $\lambda$ is plugged back into $T(\cdot)$ and some linear eigensolver is used to compute the eigenvector associated to the smallest eigenvalue of $T(\lambda)$. The algorithm for {\bf pmCIMa} is as follows.

\vskip 0.2cm
\textbf{pmCIMa} 
\begin{itemize}
\item[-] Given a series of disks $\Omega_i$, with center $z_i$ and radius $r$, covering $\Omega$.
\item[1.] Compute eigenvalues in $\Omega$ 
    \begin{itemize}
    \item[1.a.] $L = {\rm ceil}(\log_2(r/tol_{eps}))$.
    \item[1.b.] For $l = 1, \ldots, L$, 
    \begin{itemize}
    	\item[-] Compute the indicators of the all regions of level $l$ in parallel.
	\item[-] Uniformly divide the regions for which the indicators are larger than $tol_{ind}$ into smaller regions.
    \end{itemize}
     \item[1.c.] Use the centers of the regions as approximate eigenvalues $\lambda_i$'s.
    \end{itemize}
\item[2.] Compute the eigenvectors associated to the smallest eigenvalues $\lambda^0_i$'s of $T(\lambda_i)$'s in parallel.
\end{itemize}

\begin{remark}
When there exist eigenvalues outside $\Omega_i$ but close to it, the indicator is large. Hence, when the algorithm zooms in around an eigenvalue, there can be several regions close to the eigenvalue having larger indicators and a merge is implemented.
\end{remark}
\begin{remark}
Although it is out of the scope of the current paper, we note that, when a quadrature point is close to an eigenvalue, the linear system can be ill-conditioned. 
\end{remark}

\subsection{pmCIMb}
{\bf pmCIMb} combines SIM and Beyn's method.
In Step 1, we cover $\Omega$ with subregions $\{\Omega_i\}_{i=1}^I$ and compute $I_{\Omega_i}$
to determine if $\Omega_i$ contains eigenvalues. If $I_{\Omega_i} > tol_{ind}$, then $\Omega_i$ is saved for Step 2. Otherwise, $\Omega_i$ is discarded. 



\begin{remark}
Ideally each $\Omega_i$ is small such that it contains a few eigenvalues. This is clearly problem dependent. 
\end{remark}

In Step 2, for a disk $\Omega_i$, Beyn's method is used to compute candidate eigenvalues (and the associated eigenvectors). The approximations $C_{0, N}$ and $C_{1, N}$ for $C_0$ and $C_1$ are given respectively by
\begin{align}\label{integrals11}
C_0\approx C_{0,N}&=\frac{1}{iN}\sum_{j=1}^{N}ire^{i\theta_j} \left[ T(z+re^{i\theta_j})\backslash {V}\right], 
\\ \label{integrals21} 
C_1\approx C_{1,N}&=\frac{1}{iN}\sum_{j=1}^{N}ire^{i\theta_j} (z+r\theta_j) \left[ T(z+re^{i\theta_j})\backslash{V}\right]. 
\end{align}

Since the number of eigenvalues inside $\Omega_i$ is unknown, a rank test is used in Beyn's method. In contract, we employ a simpler criterion of using a threshold value $tol\_svd$ for the singular values. In fact, $tol\_svd$ is relevant to the locations of eigenvalues relative to $\partial \Omega_i$ and the number of quadrature points. Numerical results indicate that $tol\_svd = 10^{-6}$ is a reasonable choice for $N\ge64$ in \eqref{integrals11} and \eqref{integrals21}.

In Step 3, for a value $\lambda_i$ from Step 2, we compute the smallest eigenvalue $\lambda_{i}^0$ of $T(\lambda_i)$. If $\lambda_{i}^0 < 10^{-6}$, $\lambda_i$ is taken as an eigenvalue of $T(\cdot)$. Otherwise, $\lambda_i$ is discarded. Typical eigensovlers such as Arnodi methods can be used to compute the smallest eigenvalue of $T(\lambda_i)$.

The algorithm for {\bf pmCIMb} is as follows.

\vskip 0.2cm
\textbf{pmCIMb} 
\begin{itemize}
\item[-] Given a series of disks $\Omega_i$ with center $z_i$ and radius $r_i$, covering $\Omega$.
\item[1.] Decide which $\Omega_i$'s contain eigenvalues.
    \begin{itemize}
    \item[1.a.] Generate a random vector ${\boldsymbol f}$ with unit norm.
    \item[1.b.] Compute the indicators $I_{\Omega_i}$ in parallel. 
    \item[1.c.] Keep $\Omega_i$'s such that $I_{\Omega_i} > tol_{ind}$.
    \end{itemize}
\item[2.] Compute the candidate eigenvalues (and eigenvectors) in $\Omega_i$ in parallel.
    \begin{itemize}
    \item[2.a.] Generate a random matrix $\tilde{V}\in \mathbb C^{n,\ell}$ with $\ell$ large enough.
    \item[2.b.] Calculate $C_{0,N}$ and $C_{1,N}$.
    \item[2.c.]  Compute the singular value decomposition
            $$C_{0,N}=V\Sigma W^H,$$
        where $V\in\mathbb{C}^{n,\ell}$, $\Sigma=diag(\sigma_1,\sigma_2,\cdots,\sigma_{\ell})$, $W\in\mathbb{C}^{\ell,\ell}$.
    \item[2.d] Find $p, 0< p\leq \ell,$ such that
$$\sigma_1\geq \sigma_2\geq\cdots \geq\sigma_p>tol_{svd}>\sigma_{p+1}\approx\cdots\approx\sigma_{\ell}\approx0.$$
Take the first $p$ columns of the matrix $V$ denoted by $V_0=V(:,1:p)$. Similarly, $W_0=W(:,1:p)$, and $\Sigma_0=diag(\sigma_1,\sigma_2,\cdots,\sigma_p)$.
    \item[2.e] Compute the eigenvalues $\lambda_i$'s and eigenvectors $v_i$'s of 
    \[
    D=V_0^HC_{1,N}W_0\Sigma_0^{-1}\in\mathbb{C}^{p,p}.
    \]
    \item[2.f]  Delete those eigenvalues which are outside the square.
    \end{itemize}
\item[3.] Verification of $\lambda_i$'s in parallel.
\begin{itemize}
    \item[3.a] Compute the smallest eigenvalue $\lambda^0_i$ of $T(\lambda_i)$.
    \item[3.b] Output $\lambda_i$'s as eigenvalues of $T(\cdot)$ if $|\lambda^0_i| < tol_{eps}$. 
\end{itemize}
\end{itemize}

\begin{remark}
If there are many cores available, a deeper level of parallelism can be implemented, i.e., solve the linear systems at all quadrature points in parallel.
\end{remark}

For robustness, the choice of various threshold values intends not to miss any eigenvalues and  use Step 2.f and Step 3 to eliminate those which are not qualified.

We devote the rest of this section to the comparison of the two methods and the discussion of the choices of various parameters and some implementation details.
\begin{itemize}
\item {\bf pmCIMa} is simpler than {\bf pmCIMb} but need to solve more linear systems.
\item When dividing $\Omega$, avoid the real and imaginary axes as many problems has real or pure imaginary eigenvalues.
\item The threshold values $tol_{ind}, tol_{eps}, tol_{svd}$ are problem-dependent. They also depend on the quadrature rule, i.e. $N$, and the minimal distance of the eigenvalues to the contour. In general, a larger $N$ is better but $N=16$ seems to be enough for SIM and $N=64$ for Beyn's method. Note that the minimal distance is unknown. 
\item The value $\ell$ is problem-dependent but unknown.  The default value for $\ell$ is $20$.
\item If Step 1 of {\bf pmCIMb} indicates that a region contains eigenvalues but the algorithm computes none, use smaller subregions to cover $\Omega$ and/or increase $\ell$.
\end{itemize}
The main challenge is to balance the robustness and efficiency. Prior knowledge of the problem is always helpful. The user might want to test several different parameters and chose them accordingly for better performance.


\section{Numerical Examples}\label{NE}
We compute two nonlinear eigenvalue problems using  the two methods. The first one is a quadratic eigenvalue problem and the second one is the numerical approximation of scattering poles.
The computation is done using MATLAB R2021a on a Mac Studio (2022) with 128GB memory and an Apple M1 Ultra chip (20 cores with 16 performance and 4 efficiency).

\subsection{Quadratic Eigenvalue Problem}
Consider a quadratic eigenvalue problem 
\[
T(z) = T_0+z T_1 + z^2 T_2, \quad T_j \in \mathbb R^{100, 100}, j=0, 1, 2, 
\]
where $T_j$'s are generated by MATLAB {\it rand}. Let $\Omega = [-0.5, 0.5] \times [-0.5, 0.5]$. We uniformly divide $\Omega$ into $225$ squares and use the circumscribing disks as the input for {\bf pmCIMa} and {\bf pmCIMb}.

\begin{figure}[ht]
\begin{center}
\begin{tabular}{cc}
\resizebox{0.52\textwidth}{!}{\includegraphics{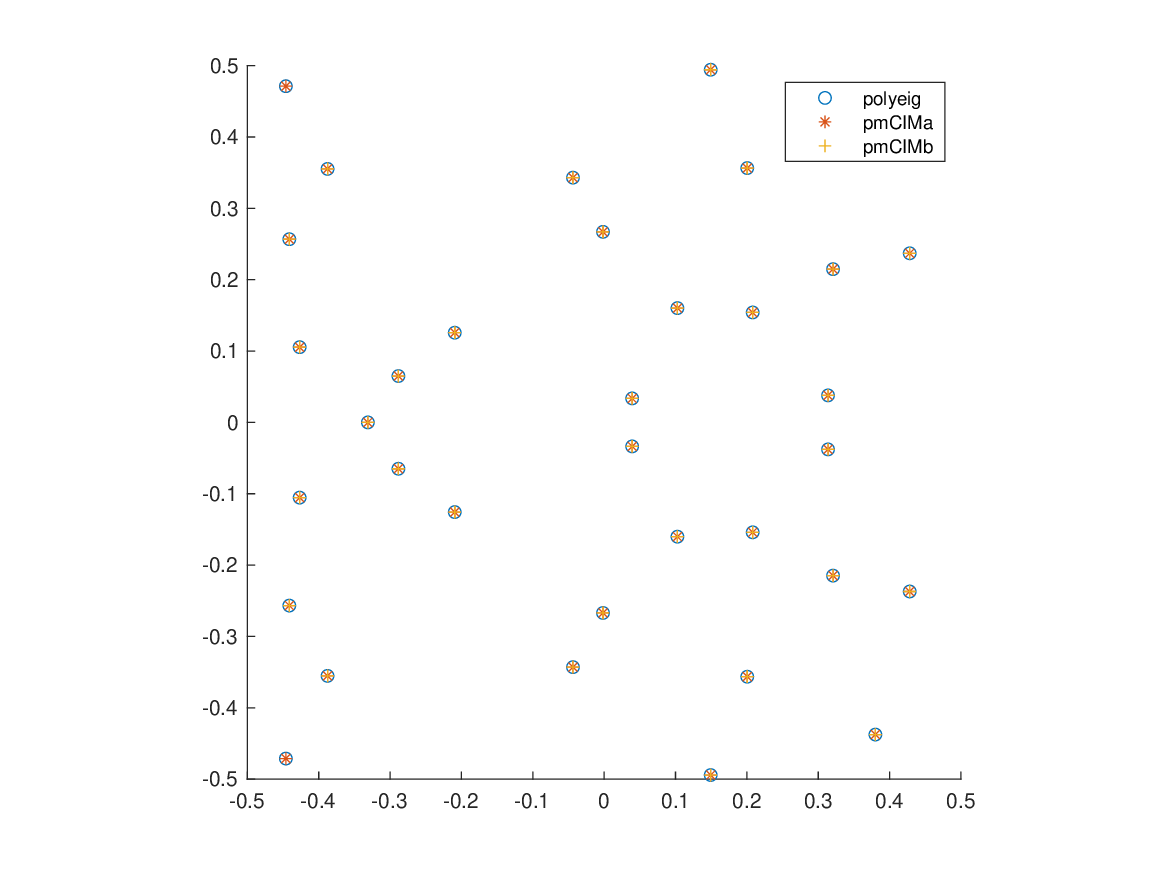}}&
\resizebox{0.52\textwidth}{!}{\includegraphics{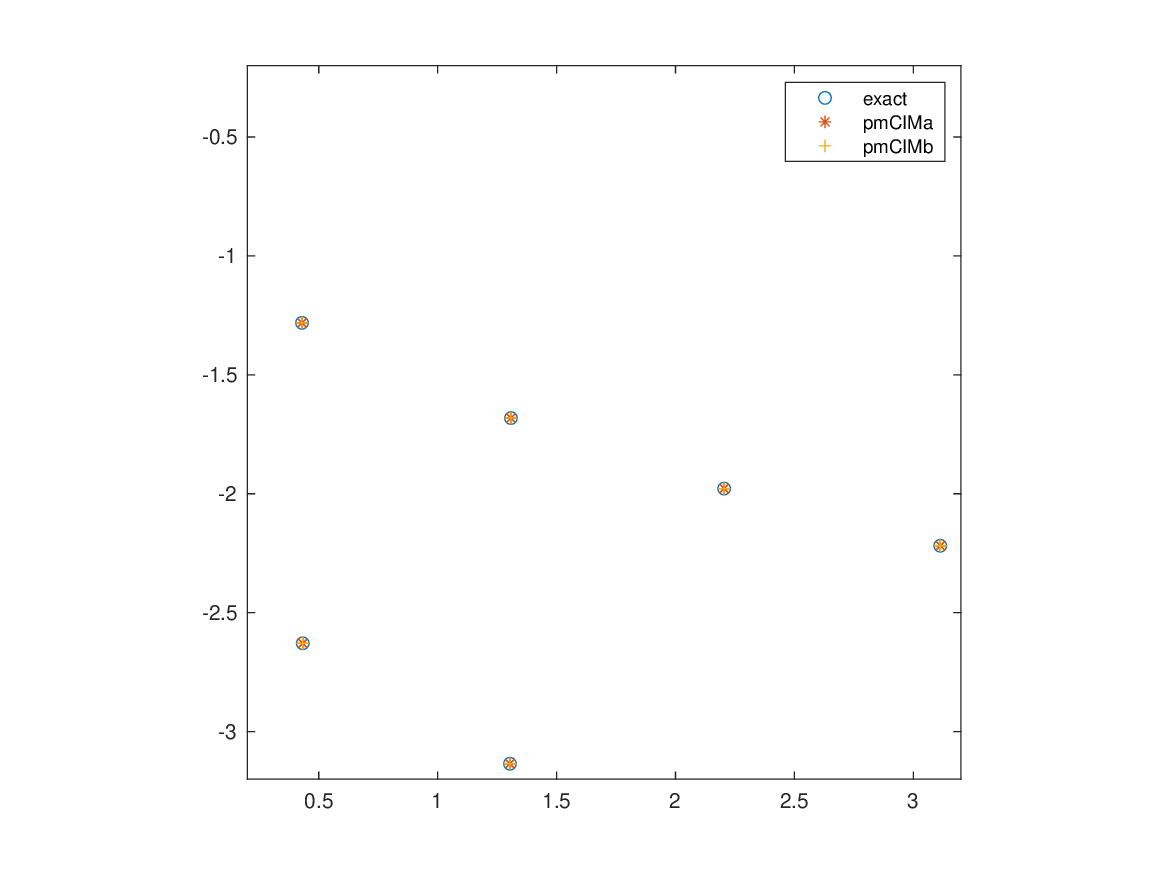}}
\end{tabular}
\end{center}
\caption{Eigenvalues computed by {\bf pmCIMa} and {\bf pmCIMb}. Left: Quadratic eigenvalue problem ($\Omega = [-0.5, 0.5] \times [-0.5, 0.5]$). Right: Scattering poles ($\Omega = [0.2, 3.2] \times [-3.2, -0.2]$). }
\label{Eigplot}
\end{figure}

The eigenvalues computed by 
MATLAB {\it polyeig}, {\bf pmCIMa} and {\bf pmCIMb} are shown in Fig.~\ref{Eigplot} (left). It can be seen that they match well. In Table~\ref{times}, we show the CPU time (in seconds) of  {\bf pmCIMa} and {\bf pmCIMb}, and their sequential versions (by replacing {\it parfor} by {\it for} in the MATLAB codes in Appendix). The parallel pool use $12$ workers and speedup is approximately $11.8$ for {\bf pmCIMa} and $10.9$ for {\bf pmCIMb}, respectively. 

\subsection{Scattering Poles}
We consider the computation of scattering poles for a sound soft unit disk $D$ \cite{MaSun2023AML}. The scattering problem is to  find $u \in H_{loc}^1(\mathbb R^2 \setminus \overline{D})$ such that
\begin{equation}\label{SbyOb}
\begin{array}{ll}
\Delta u + \kappa ^2 u = 0 & \text{in } \mathbb R^2 \setminus \overline{D}, \\
u = g & \text{on } \partial D, \\
\lim_{r \to \infty} \sqrt{r} \left(\frac{\partial u}{\partial r} - i\kappa u\right) = 0,&
\end{array}
\end{equation}
where $\kappa \in \mathbb C$ is the wave number, $g \in L^2(\partial D)$, and $r = |x|$. The scattering operator ${\mathcal B}(\kappa)$ is defined as the solution operator for \eqref{SbyOb}, i.e., $u = {\mathcal B}(\kappa) g$. 
It is well-known that $ {\mathcal B}(\kappa)$ is holomorphic on the upper half-plane of $\mathbb C$ and can be meromorphically continued to the lower half complex plane. The poles of $ {\mathcal B}(k)$ is the discrete set $\{\kappa \in {\mathbb C}: \text{Im}(\kappa) < 0\}$, which are called the scattering poles \cite{DyatlovZworski2019}. 

For $\phi \in L^2(\partial D)$, define the double layer operator $K$ as
\begin{equation}\label{DoubleS}
(K(\kappa)\phi)(x) := 2 \int_{\partial D} \frac{\partial \Phi(x, y, \kappa )}{\partial \nu(y)} \phi(y) d s(y), \quad x \in \partial D,
\end{equation}
where $\Phi(x, y, \kappa) = \frac{i}{4}H_0^{(1)}(\kappa |x-y|)$ is the fundamental solution to the Helmholtz equation and $\nu$ is the unit outward normal to $\partial D$.  Then $u$ solves \eqref{SbyOb} if it solves the integral equation
\begin{equation}\label{Integral}
\frac{1}{2}\left(I + K(\kappa) \right) \phi= g \quad \text{on } \partial D.
\end{equation}
The scattering poles are the eigenvalues of the operator $\frac{1}{2}\left(I + K(\kappa) \right)$.

We employ Nystr\"{o}m method to discretize $K(\kappa)$ to obtain an matrix $K_n(\kappa)$. Then the nonlinear eigenvalue problem is to find $\lambda$ and $\phi$ such that
\[
T(\lambda) \phi = {\boldsymbol 0} \quad \text{where} \quad T(\kappa) = \frac{1}{2}\left(I + K_n(\kappa) \right).
\]
Let $\Omega=[0.2, 3.2] \times [-3.2, -0.2]$ and divide it into $100$ uniform squares. Using the circumscribing disks as input, both {\bf pmCIMa} and {\bf pmCIMb} compute $5$ poles matching the exact poles, i.e., zeros of Hankel functions (see the right picture of Fig.~\ref{Eigplot}). In Table~\ref{times}, we show the CPU time (in seconds) of {\bf pmCIMa} and {\bf pmCIMb}, and the sequential versions (again, by replacing {\it parfor} by {\it for} in the MATLAB codes in Appendix). The parallel pool use $12$ workers and the speedup is approximately $11.1$ for {\bf pmCIMa} and $10.4$ for {\bf pmCIMb}, respectively.

\begin{table}[h!]
\label{times}
\centering
\begin{tabular}{l|r|r}
\hline
 & quadratic eigenvalue problem&scattering poles\\ \hline
mCIMa & 1098.155757 & 737.755251\\ \hline
pmCIMa  &93.184364  & 66.204724 \\ \hline
mCIMb & 77.458906 & 125.925646\\ \hline
pmCIMb  & 7.097786 & 12.138808 \\ \hline

\end{tabular}
\caption{CPU time (in seconds). {\bf mCIMa} and {\bf mCIMb} are the sequential versions for {\bf mCIMa} and {\bf mCIMb}, respectively. The speedup is almost optimal using 12 workers.}
\end{table}

\section{Conclusions}\label{CF}
We propose two parallel multi-step contour integral methods for nonlinear eigenvalue problems. Both methods use domain decomposition. Thus they are highly scalable and the speedup of the parallel algorithms is almost optimal. The effectiveness is demonstrated by a quadratic eigenvalue problem and a nonlinear eigenvalue problem of scattering poles. 

The idea of {\bf pmCIMa} is very simple and the algorithm is easy to implement. The Matlab code for {\bf pmCIMa} is less than one page.
However, {\bf pmCIMa} needs to solve more linear systems than {\bf pmCIMb} if there are a small number of eigenvalues in a large region. 
In contrast, {\bf pmCIMb} needs to have some knowledge on the number of eigenvalues inside and close to a given subregion. 
The preliminary versions of the methods have been successfully employed to compute several nonlinear eigenvalue problems of partial differential equations \cite{SunZhou2016, Xiao2021JOSAA, Gong2022MC, MaSun2023AML}.

\bibliographystyle{acmsmall}
\bibliography{mix}



\section*{Appendix - Matlab codes for pmSIM}
For questions about the codes, please send an email to jiguangs@mtu.edu.

\vskip 0.2cm
{\it test.m}. 
\begin{verbatim}
myCluster = parcluster('local');
myCluster.NumWorkers = 12; 
saveProfile(myCluster);
% compute eigenvalues in [-3,3]x[-3,3]
s = 6;
xm = -3; ym = -3; 
N = 9; 
h = s/N;
z = zeros(N^2,2);
for i= 1:N
    for j = 1:N
        z((i-1)*N+j,:)=[xm+(i-1)*h+h/2,ym+((j-1)*h+h/2)];
    end
end
r = sqrt(2)*s/N/2;
E1= pmCIMa(z,r);
E2= pmCIMb(z,r);
\end{verbatim}

{\it pmCIMa.m}
\begin{verbatim}
function allegs = pmCIMa(z,r)
% input - z,r: array of centers and radius of the disks
tol_ind = 10^(-1);      % indicator threshold
tol_eps = 10^(-6);      % eigenvalue precision
numMax = 500;           % maximum # of eigenvalues

L = ceil(log2(r/tol_eps))+2; % Level of interations
c = zeros(numMax,2);
for level=1:L
    ind = zeros(size(z,1),1);
    parfor it=1:size(z,1)
        ind(it) = indicator(z(it,:),r,16);
    end
    z0 = z(find(ind>tol_ind),:);
    if size(z0,1)>0 & level < L
        % uniformly divide the square into 4 subsquares 
        r = r/2;
        for it=1:size(z0,1)
         c((it-1)*4+1,:)=[z0(it,1)+r/sqrt(2),z0(it,2)+r/sqrt(2)];
         c((it-1)*4+2,:)=[z0(it,1)-r/sqrt(2),z0(it,2)+r/sqrt(2)];
         c((it-1)*4+3,:)=[z0(it,1)+r/sqrt(2),z0(it,2)-r/sqrt(2)];
         c((it-1)*4+4,:)=[z0(it,1)-r/sqrt(2),z0(it,2)-r/sqrt(2)];
        end
        z=c(1:4*size(z0,1),:);
    end
end
z=ones(size(z0,1),3); z(:,1:2)=z0;
finalegs=zeros(size(z0,1),2);
% merge close eigenvalues
ne=0;
for it=1:size(z0,1)
    if z(it,3)>0
       tmp=z(it,1:2);
       num=1;
       for j=1:size(z0,1)
           if (norm(z(it,1:2)-z(j,1:2))<tol_eps)
              num=num+1;
              tmp=tmp+z(j,1:2);
              z(j,3)=0;
           end
       end
       ne=ne+1;
       tmp=tmp./num;
       finalegs(ne,:)=tmp;
    end
end
allegs=finalegs(1:ne,1)+1i*finalegs(1:ne,2);
% Step 2 - validation and computation of eigenvectors
n = size(generateT(1),1);
egv = zeros(n,ne);
parfor it=1:ne
    Tn = generateT(allegs(it));
    [v,e]=eigs(Tn,1,0.001);
    egv(:,it)=v;
end
end
\end{verbatim}

{\it pmCIMa.m}
\begin{verbatim}
function finalegs = pmCIMb(z,r)
% input - z,r: array of centers and radius of the disks
tol_ind = 0.1;        % indicator threshold
tol_eps = 10^(-6);      % eigenvalue precision
tol_svd = 10^(-6);      % svd threshold  
numEig = 10;            % maximum eigenvalue in a disk
numMax = 500;           % maximum # of all eigenvalues in Omega

% Step 1 - screening 
ind = zeros(size(z,1),1);
parfor it=1:size(z,1)
    ind(it) = indicator(z(it,:),r,16);
end
z0 = z(find(ind>tol_ind),:);

size(z0,1)
% Step 2 - compute eigenvalues
egsmatrix = zeros(numEig,size(z0,1));
eignum = zeros(size(z0,1),1);
parfor it=1:size(z0,1)
    locegs = zeros(10,1);
    egs = ContourSVD(z0(it,:),r,64,numEig,tol_svd);
    locnum = length(egs);
    locegs(1:locnum) = egs;
    egsmatrix(:,it)=locegs;
    eignum(it) = locnum;
end

% exclude eigenvalues outside the square
n = 0;
allegs=zeros(numMax,1);
for it=1:size(z0,1)
    egs = egsmatrix(:,it);
    locnum = eignum(it);
    if locnum>0
        for ij=1:locnum
            if abs(real(egs(ij)) -z0(it,1))<r/sqrt(2) & abs(imag(egs(ij)) -z0(it,2))<r/sqrt(2)
                allegs(n+1)=egs(ij);
                n = n+1;
            end
        end
    end
end
allegs=allegs(1:n);

% Step 3 - validation
egsind = zeros(n,1);
parfor it=1:n
    Tn = generateT(allegs(it));
    e=eigs(Tn,1,0.001);
    if abs(e) < tol_eps
        egsind(it)=1;
    end
end
finalegs = allegs(find(egsind>0));
end

\end{verbatim}

{\it other subroutines}
\begin{verbatim}
function ind = indicator(z0,r,N)
% compute the indicator of the disc centered at z0 with radius r
M = length(generateT(1));
t=linspace(0,2*pi,N+1);
z = [r*cos(t(1:N))+z0(1)*ones(1,N);r*sin(t(1:N))+z0(2)*ones(1,N)];

f = rand(M,1); f = f/norm(f); % random vector
D0 = zeros(M,1);
D1 = zeros(M,1);
for it=1:N
    k=z(1,it)+sqrt(-1)*z(2,it);
    Tn = generateT(k);
    d0=linearsolver(Tn,f);
    D0=D0+d0*exp(1i*t(it))*r/N;
    if mod(it,2)==0
        D1=D1+2*d0*exp(1i*t(it))*r/N;
    end
end
ind = norm(D0./D1)/sqrt(M);
end

function Eigs = ContourSVD(z0,r,N,numEig,tol_svd)
M = length(generateT(1));
t=linspace(0,2*pi,N+1);
z = [r*cos(t(1:N))+z0(1)*ones(1,N);r*sin(t(1:N))+z0(2)*ones(1,N)];

f = rand(M, numEig);
for it=1:numEig
    rv = rand(M,1);
    f(:,it) = rv/norm(rv);
end
D0=zeros(M, numEig);
D1=zeros(M, numEig);
for it=1:N
    k=z(1,it)+sqrt(-1)*z(2,it);
    Tn = generateT(k);
    d0=linearsolver(Tn,f);
    d1=k*d0;
    D0=D0+d0*exp(1i*t(it))*r/N/1i;
    D1=D1+d1*exp(1i*t(it))*r/N/1i;
end

[V,E,W] = svd(D0);
sigmma=diag(E);
ind=find(sigmma>tol_svd);
p=length(ind);
V0=V(1:M,1:p);
W0=W(1:numEig,1:p);
E0=1./sigmma(1:p);
B=V0'*D1*W0*diag(E0);
Eigs = eig(B);
end

function x = linearsolver(Tn,f)
    x = Tn\f;
end

function T = generateT(z)
% Replace it with a user-defined function to generate T(z) 
 T2=diag([3 1 3 1]); 
 T1=[0.4 0 -0.3 0; 0 0 0 0;-0.3 0 0.5 -0.2;0 0 -0.2 0.2];
 T0=[-7 2 4 0; 2 -4 2 0;4 2 -9 3; 0 0 3 -3];
 T = T0+z*T1+z^2*T2;
end

\end{verbatim}
\end{document}